\RequirePackage{fix-cm}
\documentclass[smallextended]{svjour3}       
\smartqed  
\usepackage{graphicx}
\begin{document}

\title{On weak-strong uniqueness property for the full compressible magnetohydrodynamics flows 
}


\author{Weiping Yan   
}


\institute{College of of Mathematics, Jilin University, Changchun
130012, P.R. China.\\
Beijing International Center for Mathematical Research, Peking University, Beijing 100871, P.R. China.\\
\email{yan8441@126.com}}

\date{Received: date / Accepted: date}

\maketitle

\begin{abstract}
This paper is devoted to the study of the weak-strong uniqueness property for the full compressible magnetohydrodynamics flows.
The governing equations for magnetohydrodynamic flows are expressed by the
full Navier-Stokes system for compressible fluids enhanced by
forces due to the presence of the magnetic field as well as
the gravity and with an additional equation which describes the
evolution of the magnetic field. Using the relative entropy inequality, we prove that a weak solution coincides with the strong solution, emanating from the same initial data, as long as the latter exists.

\end{abstract}

\section{Introduction and Main results}
\label{Sec:1} \indent
This paper studies the weak-strong uniqueness property of the viscous compressible magnetohydrodynamic flows
\begin{eqnarray}\label{E1-1}
&&\partial_t\rho+\textbf{div}_x(\rho\textbf{u})=0,\\
\label{E1-2}
&&\partial_t(\rho\textbf{u})+\textbf{div}_x(\rho\textbf{u}\otimes\textbf{u})+\nabla_xP(\rho,\theta)=\textbf{div}_x\textbf{S}+\textbf{J}\times\textbf{H},\\
\label{E1-1'}
&&\partial_t(\rho s(\rho,\theta))+\textbf{div}_x(\rho s(\rho,\theta)\textbf{u})+\textbf{div}_x(\frac{\textbf{q}}{\theta})=\sigma,\\
\label{E1-1R}
&&\partial_t\textbf{H}-\nabla\times(\textbf{u}\times\textbf{H})+\nabla\times(\nu\nabla\times\textbf{H})=0.
\end{eqnarray}
where $\textbf{u}$ is the vector field, $\rho$ is the density, $\theta$ is the temperature, $\textbf{J}$ is the electronic current, $e(\rho,\theta)$ is the (specific) internal energy and $\textbf{H}$ is
the magnetic field. The electronic
current satisfies Amp\`{e}re's law
\begin{eqnarray}\label{E1-3}
\textbf{J}=\nabla\times\textbf{H},
\end{eqnarray}
whereas the Lorentz force is given by
\begin{eqnarray}\label{E1-4}
\textbf{J}\times\textbf{H}=\textbf{div}_x(\frac{1}{\mu}\textbf{H}\otimes\textbf{H}-\frac{1}{2\mu}|\textbf{H}|^2\textbf{I}),
\end{eqnarray}
with $\mu$ being a permeability constant of free space, which here is assumed to be $\mu=1$ for simplicity
of the presentation. The electronic current $\textbf{J}$, the electric field $\textbf{E}$ and the magnetic field $\textbf{H}$ are related
through Ohm's law
\begin{eqnarray}\label{E1-5}
\textbf{J}=\sigma(\textbf{E}+\textbf{u}\times\textbf{H}).
\end{eqnarray}
The interaction described by the theory of magnetohydrodynamics, ``collective effects,'' is governed by
the Faraday's law,
\begin{eqnarray}\label{E1-6}
\partial_t\textbf{H}+\nabla\times\textbf{E}=0,~~\textbf{div}_x\textbf{H}=0.
\end{eqnarray}
Taking into consideration (\ref{E1-5}) we are able to write (\ref{E1-6}) in the following form
\begin{eqnarray}\label{E1-7}
\partial_t\textbf{H}+\nabla\times(\textbf{H}\times\textbf{u})+\nabla\times(\nu\nabla\times\textbf{H})=0,
\end{eqnarray}
where $\nu=\frac{1}{\sigma}$.

Motivated by several recent studies devoted to the scale analysis as well as numerical experiments
related to the proposed model (see Klein et al. \cite{Klein}), we suppose that the viscous stress $\textbf{S}$ is a linear function of the velocity gradient, therefore described by Newton's law
\begin{eqnarray}\label{E1-8}
\textbf{S}(\theta,\nabla_x\textbf{u})=\mu(\theta)(\nabla_x\textbf{u}+\nabla_x^{\perp}\textbf{u}-\frac{2}{3}\textbf{div}_x\textbf{u}\textbf{I})+\eta(\theta)\textbf{div}_x\textbf{u}\textbf{I},
\end{eqnarray}
while $\textbf{q}$ is the heat flux satisfying Fourier's law
\begin{eqnarray}\label{E1-9}
\textbf{q}=-\kappa(\theta)\nabla_x\theta,
\end{eqnarray}
and $\sigma$ stands for the entropy production rate which is non-negative measure given by
\begin{eqnarray}\label{E1-10}
\sigma\geq\frac{1}{\theta}(\textbf{S}(\theta,\nabla_x\textbf{u}):\nabla_x\textbf{u}-\frac{\textbf{q}(\theta,\nabla_x\theta)\cdot\nabla_x\theta}{\theta}).
\end{eqnarray}
We supplement compressible magnetohydrodynamic flows (\ref{E1-1})-(\ref{E1-1R}) with conservation boundary condition:
\begin{eqnarray}\label{E1-11}
\textbf{u}|_{\partial\Omega}=\textbf{q}\cdot\textbf{n}|_{\partial\Omega}=0,
\end{eqnarray}
and
\begin{eqnarray}\label{E1-12}
\textbf{H}|_{\partial\Omega}=0.
\end{eqnarray}
The concept of weak solution in fluid dynamics was introduced by Leray \cite{Ler} in the context of incompressible, linearly viscous fluids.
The original ideas of Leray have been put into the elegant framework of generalized derivatives
(distributions) and the associated abstract function spaces of Sobolev type (For example, see Ladyzhenskaya \cite{La} and Temam \cite{T}).
Lions \cite{Lions} extended the theory to the class of barotropic flows (see also \cite{Fei0}). One of meaningful compressible flow models is the compressible magnetohydrodynamics (MHD). It is a combination of the compressible Navier-Stokes equations of fluid dynamics and Maxwell's equations of electromagnetism.
Ducomet and Feireisl \cite{Du} proved that the existence of global in time weak solutions to
a multi-dimensional nonisentropic MHD system for gaseous stars coupled with the Poisson equation with
all the viscosity coefficients and the pressure depending on temperature and density asymptotically, respectively.
Hu and Wang \cite{Hu1} studied the global variational weak
solution to the three-dimensional full magnetohydrodynamic equations
with large data by an approximation scheme and a weak convergence
method. Jiang, et all. \cite{Jiang} obtained that the convergence towards the strong solution of the
ideal incompressible MHD system in the periodic domains. Recently, Kwon, et all \cite{Kwon} established the incompressible limits of weak solutions to the compressible magnetohydrodynamics flows (\ref{E1-1})-(\ref{E1-1R}) on both
bounded and unbounded domains.

The physical properties of the magnetohydrodynamics flows are reflected through various constitutive
relations which are expressed as typically non-linear functions relating the pressure $P=P(\rho,\theta)$,
the internal energy $e(\rho,\theta)$, the specific entropy $s=s(\rho,\theta)$ to the macroscopic variables $\rho$, $\textbf{u}$,
and $\theta$. According to the fundamental principles of thermodynamics, the specific internal energy $e$ is
related to the pressure $P$, and the specific entropy s through Gibbs' relation
\begin{eqnarray}\label{E1-13}
\theta Ds(\rho,\theta)=De(\rho,\theta)+P(\rho,\theta)D(\frac{1}{\rho}),
\end{eqnarray}
where $D$ denotes the differential with respect to the state variables $\rho$ and $\theta$.

Since the lack of information resulting from the inequality sign in (\ref{E1-10}), we need supplement the resulting system with the energy inequality,
\begin{eqnarray}\label{E1-14}
\frac{d}{dt}\int_{\Omega}(\frac{1}{2}\rho|\textbf{u}|^2+\rho e(\rho,\theta)+\frac{1}{2}|\textbf{H}|^2)dx+\int_{\Omega}(|\nabla\textbf{u}|^2+\nu|\nabla\times\textbf{H}|^2)\leq0.
\end{eqnarray}
Thus the total energy $\mathcal{E}$ is given by
\begin{eqnarray}\label{E1-14R}
\mathcal{E}=\frac{1}{2}\rho|\textbf{u}|^2+\rho e(\rho,\theta)+\frac{1}{2}|\textbf{H}|^2.
\end{eqnarray}
Under these circumstances, it
can be shown (see \cite{Fei2}, Chapter 2) that any weak solution of (\ref{E1-1}) that is sufficiently
smooth satisfies, instead of (\ref{E1-10}), the standard relation
\begin{eqnarray}\label{E1-10R1}
\sigma=\frac{1}{\theta}\left(\textbf{S}(\theta,\nabla_x\textbf{u}):\nabla_x\textbf{u}-\frac{\textbf{q}(\theta,\nabla_x\theta)\cdot\nabla_x\theta}{\theta}\right).
\end{eqnarray}
The pressure $P=P(\rho,\theta)$ is here expressed as
\begin{eqnarray}\label{E1-15}
P=P_F+P_R,~~P_R=\frac{a}{3}\theta^4,~a>0,
\end{eqnarray}
where $P_R$ denotes the radiation pressure. Moreover, we shall assume that $P_F=P_M+P_E$ , where $P_M$ is the classical molecular pressure
obeying Boyle's law, while $P_E$ is the pressure of electron gas constituent behaving like a Fermi
gas in the degenerate regime of high densities and/or low temperatures (see Chapters 1, 15 in
Eliezer et al. \cite{E}). Thus necessarily $P_F$ takes the form
\begin{eqnarray}\label{E1-16}
P_F=\theta^{\frac{5}{2}}p(\frac{\rho}{\theta^{\frac{3}{2}}}),
\end{eqnarray}
where $p\in\textbf{C}^1[0,\infty)$ satisfies
\begin{eqnarray}\label{E1-17}
p(0)=0,~~p'(Z)>0~for~all~z\geq0.
\end{eqnarray}
In agreement with Gibbs' relation (\ref{E1-13}), the internal energy can be taken as
\begin{eqnarray}\label{E1-18}
e=e_F+e_R,~~with~e_R=a\frac{\theta^4}{\rho},
\end{eqnarray}
where $e_F=e_F(\rho,\theta)$, $P_F(\rho,\theta)$ are interrelated through the following equation of state
\begin{eqnarray}\label{E1-18RR}
P_F(\rho,\theta)=\frac{2}{3}\rho e_F(\rho,\theta).
\end{eqnarray}
We need the thermodynamic stability hypothesis:
\begin{eqnarray}\label{E1-18R}
\frac{\partial P(\rho,\theta)}{\partial\rho}>0,~~\frac{\partial e(\rho,\theta)}{\partial\theta}>0~~for~all~~\rho,\theta>0.
\end{eqnarray}
The second inequality in thermodynamic stability hypothesis (\ref{E1-18R}) gives that
\begin{eqnarray}\label{E1-19}
0<\frac{\frac{5}{3}p(Z)-p'(Z)Z}{Z}<c~~for~all~~Z>0,
\end{eqnarray}
which implies that the function $Z\mapsto\frac{p(Z)}{Z^{\frac{5}{3}}}$ is decreasing and we suppose that
\begin{eqnarray}\label{E1-20}
\lim_{Z\longrightarrow\infty}\frac{p(Z)}{Z^{\frac{5}{3}}}=p_{\infty}>0.
\end{eqnarray}
In accordance with (\ref{E1-13}) and (\ref{E1-18}) we set the entropy as
\begin{eqnarray}\label{E1-21}
s=s_F+s_R,~~wtih~s_F=S(\frac{\rho}{\theta^{\frac{3}{2}}}),~~s_R=\frac{4a}{3\rho}\theta^3.
\end{eqnarray}
Furthermore, by the Third law of thermodynamics,
\begin{eqnarray}\label{E1-22}
S'(Z)=-\frac{3}{2}\frac{\frac{5}{3}p(Z)-p'(Z)Z}{Z^2}<0,~~\lim_{Z\longrightarrow\infty}S(Z)=0.
\end{eqnarray}
We choose the transport coefficients in the form
\begin{eqnarray}\label{E1-23}
&&\mu(\theta)=\mu_0+\mu_1\theta,~~\mu_0,\mu_1>0,~~\eta\equiv0,\\
\label{E1-24}
&&\kappa(\theta)=\kappa_0+\kappa_2\theta^2+\kappa_3\theta^3,~~\kappa_i>0,~~i=0,2,3.
\end{eqnarray}

A fundamental test of admissibility of a class of weak solutions to a given evolutionary
problem is the property of weak-strong uniqueness. More specifically, the
weak solution must coincide with a (hypothetical) strong solution emanating from the same initial data as long as the latter exists. This problem has been intensively
studied for the incompressible Navier-Stokes system, for example, see \cite{Che,Pi1,Ser}. It is a bit more delicate in the case of compressible cases. The weak-
strong uniqueness of compressible barotropic Navier-Stokes system and isentropic compressible Navier-Stokes system were established in \cite{Fei11,Fei3} and \cite{Pi2}, respectively. P. Germain \cite{Pi2} provides only a partial
and conditional answer to the weak-strong uniqueness problem for the compressible Navier-Stokes equations. This question is definitely solved in \cite{Fei3}.
More recently, Feireisl and Novotn\'{y} \cite{Fei1} extended the problem to compressible Navier-Stokes-Fourier system by the relative entropy inequality.
The relative entropy in \cite{Fei1} is reminiscent to C.M. Dafermos \cite{Daf} (who introduced
the relatives entropies via the entropy flux pairs for the conservation laws),
but is different from the C.M. Dafermos concept (in contrast to \cite{Daf}, it is based on
the thermodynamic stability conditions).

Inspired by the work of Feireisl and Novotn\'{y} \cite{Fei1},
we prove that the weak-
strong uniqueness of compressible three-dimensional magnetohydrodynamic equations. Our contribution is to
construct suitable relative entropy inequality to (\ref{E1-1})-(\ref{E1-1'}). Then we overcome the presence of the magnetic field and its interaction with the
hydrodynamic motion in the MHD flow of large oscillation.

We organize the rest of this paper as follows. In section 2, we recall the definition of
the weak solutions and strong solutions to the magnetohydrodynamic
flows on bounded domains. Meanwhile, the relative entropy inequality of (\ref{E1-1})-(\ref{E1-1R}) is derived.
In the last section, we give the rigorous proof of the weak-strong uniqueness property for the compressible magnetohydrodynamic
flows on bounded domains in the spirit of Feireisl and Novotn\'{y} \cite{Fei1}.

\section{Relative entropy and Main result}
Let $\Omega\subset\textbf{R}^3$ be a bounded Lipschitz domain. We recall the definition of weak solution for (\ref{E1-1})-(\ref{E1-1R}).
\begin{definition}
We say that a quantity $(\rho,\textbf{u},\theta,\textbf{H})$ is a weak solution of the
full magnetohydrodynamic flows (MHD) (\ref{E1-1})-(\ref{E1-1R}) supplemented with the initial data $(\rho_0,\textbf{u}_0,s(\rho_0,\theta_0),\textbf{H}_0)$, and $\rho_0\geq0,\theta_0>0$ provided
that the following holds.

\textbf{i)} The density $\rho$ is a non-negative function, $\rho\in\textbf{C}_{weak}([0,T];\textbf{L}^{\frac{5}{3}}(\Omega))$, the velocity field $\textbf{u}\in\textbf{L}^2(0,T;\textbf{W}_0^{1,2}(\Omega;\textbf{R}^3))$, $\rho\textbf{u}\in\textbf{C}_{weak}([0,T];\textbf{L}^{\frac{5}{4}}(\Omega;\textbf{R}^3))$. Equation (\ref{E1-1}) is replaced by a family of integral identities
\begin{eqnarray}\label{E2-1}
\int_{\Omega}\rho(\tau,\cdot)\varphi(\tau,\cdot)dx-\int_{\Omega}\rho_0\varphi(0,\cdot)dx=\int_0^{\tau}\int_{\Omega}(\rho\partial_t\varphi+\rho\textbf{u}\cdot\nabla_x\varphi)dxdt
\end{eqnarray}
for any $\varphi\in\textbf{C}^1([0,T]\times\bar{\Omega})$, and any $\tau\in[0,T]$.

\textbf{ii)} The balance of momentum holds in distributional sense, namely
\begin{eqnarray}\label{E2-2}
&&\int_{\Omega}\rho\textbf{u}(\tau,\cdot)\cdot\varphi(\tau,\cdot)dx-\int_{\Omega}\rho_0\textbf{u}_0\cdot\varphi(0,\cdot)dx\nonumber\\
&=&\int_0^{\tau}\int_{\Omega}(\rho\textbf{u}\cdot\partial_t\varphi+\rho\textbf{u}\otimes\textbf{u}:\nabla_x\varphi+P\textbf{div}_x\varphi-\textbf{S}:\nabla_x\varphi\nonumber\\
&&+[(\nabla\times\textbf{H})\times\textbf{H}]\cdot\varphi)dxdt
\end{eqnarray}
for any $\varphi\in\textbf{C}^1([0,T]\times\bar{\Omega};\textbf{R}^3)$, $\varphi|_{\partial\Omega}=0$ and any $\tau\in[0,T]$.

\textbf{iii)} The entropy balance (\ref{E1-1'}) and (\ref{E1-10}) are replaced by a family of integral inequalities
\begin{eqnarray}\label{E2-3}
&&\int_{\Omega}\rho s(\rho_0,\theta_0)\varphi(0,\cdot)dx-\int_{\Omega}\rho s(\rho,\theta)(\tau,\cdot)\varphi(\tau,\cdot)dx\nonumber\\
&+&\int_0^{\tau}\int_{\Omega}\left(\frac{\varphi}{\theta}(\textbf{S}:\nabla_x\textbf{u}-\frac{\textbf{q}\cdot\nabla_x\theta}{\theta}\right)dxdt\nonumber\\
&\leq&-\int_0^{\tau}\int_{\Omega}\left(\rho s(\rho,\theta)\partial_t\varphi+\rho s(\rho,\theta)\textbf{u}\cdot\nabla_x\varphi+\frac{\textbf{q}\cdot\nabla_x\varphi}{\theta}\right)dxdt
\end{eqnarray}
for any $\varphi\in\textbf{C}^1([0,T]\times\bar{\Omega})$, $\varphi\geq0$ and almost all $\tau\in[0,T]$. Here the quantities $\textbf{S}$ and $\textbf{q}$ are given through the constitutive equations (\ref{E1-8}) and (\ref{E1-9}). Moreover,
similarly to the above, all quantities must be at least integrable on $(0,T)\times\Omega$. In particular,
$\theta$ belongs to $\textbf{L}^{\infty}(0,T;\textbf{L}^4(\Omega))\cap\textbf{L}^2(0,T;\textbf{W}^{1,2}(\Omega))$. In addition, we require $\theta$ to be positive for almost all $(t,x)\in(0,T)\times\Omega$.

\textbf{iv)} The total energy of the system satisfies the following inequality
\begin{eqnarray}\label{E2-4}
\int_{\Omega}\left(\frac{1}{2}\rho|\textbf{u}|^2+\rho e(\rho,\theta)+\frac{1}{2}|\textbf{H}|^2\right)dx
&+&\int_0^{\tau}\int_{\Omega}(|\nabla\textbf{u}|^2+\nu|\nabla\times\textbf{H}|^2)dxdt\nonumber\\
&\leq&\int_{\Omega}\left(\frac{1}{2}\rho_0|\textbf{u}_0|^2+\rho_0e(\rho_0,\theta_0)+\frac{1}{2}|\textbf{H}_0|^2\right)dx~~~~~~~~
\end{eqnarray}
for almost all $\tau\in[0,T]$.

\textbf{v)} The magnetic field $\textbf{H}\in\textbf{L}^2(0,T;\textbf{W}^{1,2}(\Omega;\textbf{R}^3))$. The Maxwell equation (\ref{E1-1R}) verifies
\begin{eqnarray}\label{E2-5}
\int_{\Omega}\textbf{H}(\tau,\cdot)\varphi(\tau,\cdot)dx&-&\int_{\Omega}\textbf{H}_0\varphi_0dx\nonumber\\
&=&\int_0^{\tau}\int_{\Omega}\left(\textbf{H}\cdot\partial_t\varphi-(\textbf{H}\times\textbf{u}+\nu\nabla\times\textbf{H})\cdot(\nabla\times\varphi)\right)dxdt,~~~~~~~~
\end{eqnarray}
where $\varphi\in\textbf{C}^1([0,T]\times\bar{\Omega};\textbf{R}^3)$, $\varphi|_{\partial\Omega}=0$ and any $\tau\in[0,T]$.
\end{definition}
The definition of strong solution is
\begin{definition}
We say that $(\rho',\textbf{u}',\theta',\textbf{H}')$ is a classical (strong) solution to the full magnetohydrodynamic system (\ref{E1-1})-(\ref{E1-1R}) in $(0,T)\times\Omega$ if
\begin{eqnarray}\label{E2-x1}
&&\rho'\in\textbf{C}^1([0,T]\times\bar{\Omega}),~~\theta',\partial_t\theta',\nabla^2\theta'\in\textbf{C}([0,T]\times\Omega),\nonumber\\
&&\textbf{u}',\partial_t\textbf{u}',\nabla^2\textbf{u}'\in\textbf{C}([0,T]\times\Omega;\textbf{R}^3),~~
\textbf{H}',\partial_t\textbf{H}',\nabla^2\textbf{H}'\in\textbf{C}([0,T]\times\Omega;\textbf{R}^3),~~~~~\nonumber\\
&&\rho'(t,x)\geq\rho>0,~~\theta'(t,x)\geq\theta'_0>0,~~for~all~(t,x),
\end{eqnarray}
and $\rho',\textbf{u}',\theta',\textbf{H}'$ satisfy equations (\ref{E1-1})-(\ref{E1-1R}), (\ref{E1-10R1}), together with the boundary conditions
(\ref{E1-11})-(\ref{E1-12}). Observe that hypothesis (\ref{E2-x1}) implies the following regularity properties of
the initial data:
\begin{eqnarray}\label{E2-x2}
&&\rho(0)=\rho_0\in\textbf{C}^1(\bar{\Omega}),~~\rho_0\geq\rho_0'>0,\nonumber\\
&&\textbf{u}(0)=\textbf{u}_0\in\textbf{C}^2(\bar{\Omega}),\nonumber\\
&&\theta(0)=\theta_0\in\textbf{C}^2(\bar{\Omega}),~~\theta_0\geq\theta_0'>0,\nonumber\\
&&\textbf{H}(0)=\textbf{H}_0\in\textbf{C}^2(\bar{\Omega}).
\end{eqnarray}
\end{definition}

Before giving the main result, we deduce a relative entropy inequality which is satisfied by any weak solution to the full magnetohydrodynamic system (\ref{E1-1})-(\ref{E1-1R}).

Let $\{A,B,C,D\}$ be a quantity of smooth function, $A$ and $C$ bounded below away from zero in $[0,T]\times\Omega$, and $B|_{\partial\Omega}=D|_{\partial\Omega}=0$. Moreover, we assume that smooth functions $B$ and $D$ satisfy that
\begin{eqnarray}\label{E2-5R1}
\partial_tD-\nabla\times(B\times D)+\nabla\times(\nu\nabla\times D)=0.
\end{eqnarray}
Taking $\varphi=\frac{1}{2}|B|^2$, $\varphi=B$ and $\varphi=C>0$ as a test function in (\ref{E2-1}), (\ref{E2-2}) and the entropy inequality (\ref{E2-3}), respectively, we get
\begin{eqnarray}\label{E2-x3}
\int_{\Omega}\frac{1}{2}\rho|B|^2(\tau,\cdot)dx-\int_{\Omega}\frac{1}{2}\rho_0|B|^2(0,\cdot)dx=\int_{0}^{\tau}\int_{\Omega}(\rho B\cdot\partial_t B+\rho\textbf{u}\cdot\nabla_xB\cdot B)dxdt,~~~~
\end{eqnarray}
\begin{eqnarray}\label{E2-x4}
\int_{\Omega}\rho\textbf{u}\cdot B(\tau,\cdot)dx-\int_{\Omega}\rho_0\textbf{u}_0\cdot B(0,\cdot)dx&=&\int_{0}^{\tau}\int_{\Omega}(\rho\textbf{u}\cdot\partial_tB
+\rho\textbf{u}\otimes\textbf{u}:\nabla_x B+P(\rho,\theta)\textbf{div}_xB\nonumber\\
&&-\textbf{S}(\theta,\nabla_x\textbf{u}):\nabla_xB+((\nabla\times\textbf{H})\times\textbf{H})\cdot B)dxdt~~~~~~
\end{eqnarray}
and
\begin{eqnarray}\label{E2-x5}
\int_{\Omega}\rho_0 s(\rho_0,\theta_0)C(0,\cdot)dx&-&\int_{\Omega}\rho s(\rho,\theta)C(\tau,\cdot)dx+\int_{0}^{\tau}\int_{\Omega}\frac{C}{\theta}(\textbf{S}(\theta,\nabla_x\textbf{u}):\nabla_x\textbf{u}-\frac{\textbf{q}(\theta,\nabla_x\theta)\cdot\nabla_x\theta}{\theta})dxdt\nonumber\\
&\leq&-\int_{0}^{\tau}\int_{\Omega}(\rho s(\rho,\theta)\partial_tC+\rho s(\rho,\theta)\textbf{u}\cdot\nabla_x C+\frac{\textbf{q}(\theta,\nabla_x\theta)}{\theta}\cdot\nabla_x C)dxdt.~~~~~
\end{eqnarray}
It follows from (\ref{E2-x3}), (\ref{E2-x4}) and the energy inequality (\ref{E2-4}) that
\begin{eqnarray}\label{E2-x6}
\int_{\Omega}(\frac{1}{2}\rho|\textbf{u}-B|^2&+&\rho e(\rho,\theta)+\frac{1}{2}|\textbf{H}|^2)(\tau,\cdot)dx
+\int_0^{\tau}\int_{\Omega}(|\nabla\textbf{u}|^2+\nu|\nabla\times\textbf{H}|^2)dxdt\nonumber\\
&\leq&\int_{\Omega}(\frac{1}{2}\rho_0|\textbf{u}_0-B(0,\cdot)|^2+\rho_0 e(\rho_0,\theta_0)+\frac{1}{2}|\textbf{H}_0|^2)dx\nonumber\\
&&+\int_{0}^{\tau}\int_{\Omega}((\rho\partial_tB+\rho\textbf{u}\cdot\nabla_xB)\cdot(B-\textbf{u})-P(\rho,\theta)\textbf{div}_xB\nonumber\\
&&+\textbf{S}(\theta,\nabla_x\textbf{u}):\nabla_xB-((\nabla\times\textbf{H})\times\textbf{H})\cdot B)dxdt.
\end{eqnarray}
Then summing up (\ref{E2-x5}) and (\ref{E2-x6}), we deduce that
\begin{eqnarray}\label{E2-x6'}
&&\int_{\Omega}(\frac{1}{2}\rho|\textbf{u}-B|^2+\rho e(\rho,\theta)+\frac{1}{2}|\textbf{H}|^2-C\rho s(\rho,\theta))(\tau,\cdot)dx\nonumber\\
&&+\int_0^{\tau}\int_{\Omega}(|\nabla\textbf{u}|^2+\nu|\nabla\times\textbf{H}|^2)dxdt+\int_{0}^{\tau}\int_{\Omega}\frac{C}{\theta}(\textbf{S}(\theta,\nabla_x\textbf{u}):\nabla_x\textbf{u}-\frac{\textbf{q}(\theta,\nabla_x\theta)\cdot\nabla_x\theta}{\theta})dxdt\nonumber\\
&\leq&\int_{\Omega}(\frac{1}{2}\rho_0|\textbf{u}_0-B(0,\cdot)|^2+\rho_0 e(\rho_0,\theta_0)+\frac{1}{2}|\textbf{H}_0|^2+C(0,\cdot)\rho_0 s(\rho_0,\theta_0))dx\nonumber\\
&&+\int_{0}^{\tau}\int_{\Omega}((\rho\partial_tB+\rho\textbf{u}\cdot\nabla_xB)\cdot(B-\textbf{u})-P(\rho,\theta)\textbf{div}_xB+\textbf{S}(\theta,\nabla_x\textbf{u}):\nabla_xB-((\nabla\times\textbf{H})\times\textbf{H})\cdot B)dxdt\nonumber\\
&&-\int_{0}^{\tau}\int_{\Omega}(\rho s(\rho,\theta)\partial_tC+\rho s(\rho,\theta)\textbf{u}\cdot\nabla_x C+\frac{\textbf{q}(\theta,\nabla_x\theta)}{\theta}\cdot\nabla_x C)dxdt.
\end{eqnarray}
Taking a test function $\varphi=D$ in (\ref{E2-5}) and $\varphi=\partial_{\rho}H_{C}(A,C)$ in (\ref{E2-1}), we have
\begin{eqnarray}\label{E2-5R}
\int_{\Omega}\textbf{H}(\tau,\cdot)D(\tau,\cdot)dx-\int_{\Omega}\textbf{H}_0D_0dx
=\int_0^{\tau}\int_{\Omega}\left(\textbf{H}\cdot\partial_tD-(\textbf{H}\times\textbf{u}+\nu\nabla\times\textbf{H})\cdot(\nabla\times D)\right)dxdt,~~~~
\end{eqnarray}
\begin{eqnarray}\label{E2-x7}
\int_{\Omega}\rho\partial_{\rho}H_{C}(A,C)(\tau,\cdot)dx&-&\int_{\Omega}\rho_0\partial_{\rho}H_{C(0,\cdot)}(A(0,\cdot),C(0,\cdot))dx\nonumber\\
&=&\int_0^{\tau}\int_{\Omega}(\rho\partial_t(\partial_{\rho}H_{C}(A,C)))+\rho\textbf{u}\cdot\nabla_x(\partial_{\rho}H_{C}(A,C))dxdt,~~~~~~
\end{eqnarray}
Multiplying (\ref{E2-5R1}) by $D$ and integrate over $(0,\tau)\times\Omega$, we find
\begin{eqnarray}\label{E2-5R2}
\int_{\Omega}\frac{1}{2}|D|^2(\tau,\cdot)dx-\int_{\Omega}\frac{1}{2}|D_0|^2dx
=-\int_0^{\tau}\int_{\Omega}\left((D\times B+\nu\nabla\times D)\cdot(\nabla\times D)\right)dxdt.~~~~
\end{eqnarray}
So by (\ref{E2-x6'})-(\ref{E2-5R2}), we have
\begin{eqnarray}\label{E2-x8}
&&\int_{\Omega}(\frac{1}{2}\rho|\textbf{u}-B|^2+\frac{1}{2}|\textbf{H}-D|^2+H_{C}(\rho,\theta)-\partial_{\rho}(H_{C})(A,C)(\rho-A)-H_{C}(A,C))(\tau,\cdot)dx\nonumber\\
&&+\int_{0}^{\tau}\int_{\Omega}\left(|\nabla\textbf{u}|^2+\nu(|\nabla\times D|^2-|\nabla\times D||\nabla\times\textbf{H}|+|\nabla\times\textbf{H}|^2)\right)dxdt\nonumber\\
&&+\int_{0}^{\tau}\int_{\Omega}\frac{C}{\theta}(\textbf{S}(\theta,\nabla_x\textbf{u}):\nabla_x\textbf{u}-\frac{\textbf{q}(\theta,\nabla_x\theta)\cdot\nabla_x\theta}{\theta})dxdt\nonumber\\
&\leq&\int_{\Omega}(\frac{1}{2}\rho_0|\textbf{u}_0-B(0,\cdot)|^2+\frac{1}{2}|\textbf{H}_0-D_0|^2+(H_{C(0,\cdot)}(\rho(0,\cdot),\theta(0,\cdot))\nonumber\\
&&-\partial_{\rho}(H_{C(0,\cdot)})(A(0,\cdot),C(0,\cdot))(\rho_0-A(0,\cdot))-H_{C(0,\cdot)}(A(0,\cdot),C(0,\cdot))))dx\nonumber\\
&&+\int_{0}^{\tau}\int_{\Omega}\left((\rho\partial_tB+\rho\textbf{u}\cdot\nabla_xB)\cdot(B-\textbf{u})-P(\rho,\theta)\textbf{div}_xB
+\textbf{S}(\theta,\nabla_x\textbf{u}):\nabla_xB\right)dxdt\nonumber\\
&&+\int_{0}^{\tau}\int_{\Omega}(-\textbf{H}\cdot\partial_tD-((\nabla\times\textbf{H})\times\textbf{H})\cdot B-(D\times B)\cdot(\nabla\times D)+(\textbf{H}\times\textbf{u})\cdot(\nabla\times D)) dxdt\nonumber\\
&&-\int_{0}^{\tau}\int_{\Omega}(\rho s(\rho,\theta)\partial_tC+\rho s(\rho,\theta)\textbf{u}\cdot\nabla_x C+\frac{\textbf{q}(\theta,\nabla_x\theta)}{\theta}\cdot\nabla_x C)dxdt\nonumber\\
&&-\int_0^{\tau}\int_{\Omega}(\rho\partial_t(\partial_{\rho}(H_{C})(A,C)))+\rho\textbf{u}\cdot\nabla_x(\partial_{\rho}H_{C}(A,C))dxdt\nonumber\\
&&+\int_0^{\tau}\int_{\Omega}\partial_t(A\partial_{\rho}(H_{C})(A,C)-H_{C}(A,C))dxdt.
\end{eqnarray}
Replacing $\partial_tD$ by (\ref{E2-5R1}) in (\ref{E2-x8}) to find
\begin{eqnarray}\label{E2-R5}
&&\int_{\Omega}(\frac{1}{2}\rho|\textbf{u}-B|^2+\frac{1}{2}|\textbf{H}-D|^2+H_{C}(\rho,\theta)-\partial_{\rho}(H_{C})(A,C)(\rho-A)-H_{C}(A,C))(\tau,\cdot)dx\nonumber\\
&&+\int_{0}^{\tau}\int_{\Omega}\left(|\nabla\textbf{u}|^2+\nu|\nabla\times D-\nabla\times\textbf{H}|^2\right)dxdt+\int_{0}^{\tau}\int_{\Omega}\frac{C}{\theta}(\textbf{S}(\theta,\nabla_x\textbf{u}):\nabla_x\textbf{u}-\frac{\textbf{q}(\theta,\nabla_x\theta)\cdot\nabla_x\theta}{\theta})dxdt\nonumber\\
&\leq&\int_{\Omega}(\frac{1}{2}\rho_0|\textbf{u}_0-B(0,\cdot)|^2+\frac{1}{2}|\textbf{H}_0-D_0|^2+\frac{1}{2}\rho_0|B_0|^2+(H_{C(0,\cdot)}(\rho(0,\cdot),\theta(0,\cdot))\nonumber\\
&&-\partial_{\rho}(H_{C(0,\cdot)})(A(0,\cdot),C(0,\cdot))(\rho_0-A(0,\cdot))-H_{C(0,\cdot)}(A(0,\cdot),C(0,\cdot))))dx\nonumber\\
&&+\int_{0}^{\tau}\int_{\Omega}\left((\rho\partial_tB+\rho\textbf{u}\cdot\nabla_xB)\cdot (B-\textbf{u}))-P(\rho,\theta)\textbf{div}_xB
+\textbf{S}(\theta,\nabla_x\textbf{u}):\nabla_xB\right)dxdt\nonumber\\
&&+\int_{0}^{\tau}\int_{\Omega}\left(-((\nabla\times\textbf{H})\times\textbf{H})\cdot B-(D\times B)\cdot(\nabla\times D)+(\textbf{H}\times\textbf{u})\cdot(\nabla\times D)\right) dxdt\nonumber\\
&&-\int_{0}^{\tau}\int_{\Omega}(\nabla\times(B\times D))\cdot\textbf{H}dxdt-\int_{0}^{\tau}\int_{\Omega}(\rho s(\rho,\theta)\partial_tC+\rho s(\rho,\theta)\textbf{u}\cdot\nabla_x C+\frac{\textbf{q}(\theta,\nabla_x\theta)}{\theta}\cdot\nabla_x C)dxdt\nonumber\\
&&-\int_0^{\tau}\int_{\Omega}(\rho\partial_t(\partial_{\rho}(H_{C})(A,C)))+\rho\textbf{u}\cdot\nabla_x(\partial_{\rho}H_{C}(A,C))dxdt\nonumber\\
&&+\int_0^{\tau}\int_{\Omega}\partial_t(A\partial_{\rho}(H_{C})(A,C)-H_{C}(A,C))dxdt.
\end{eqnarray}
Note that
\begin{eqnarray}\label{E9-5}
\int_{\Omega}\left((\nabla\times\textbf{H})\times\textbf{H}\right)\cdot Bdx=-\int_{\Omega}\left(\textbf{H}^{\top}B\textbf{H}+\frac{1}{2}\nabla(|\textbf{H}|^2)\cdot B\right)dx,
\end{eqnarray}
\begin{eqnarray*}
\int_{\Omega}\left(\nabla\times(B\times\textbf{H})\right)\cdot Bdx=\int_{\Omega}\left(\textbf{H}^{\top}B\textbf{H}+\frac{1}{2}\nabla(|\textbf{H}|^2)\cdot B\right)dx.
\end{eqnarray*}
So direct calculation shows that
\begin{eqnarray}\label{E2-8}
&&\int_{0}^{\tau}\int_{\Omega}\left(-((\nabla\times\textbf{H})\times\textbf{H})\cdot B-(D\times B)\cdot(\nabla\times D)+(\textbf{H}\times\textbf{u})\cdot(\nabla\times D)\right) dxdt\nonumber\\
&&-\int_{0}^{\tau}\int_{\Omega}\nabla\times(B\times D)\cdot\textbf{H}dxdt\nonumber\\
&=&\int_{0}^{\tau}\int_{\Omega}\left(-((\nabla\times\textbf{H})\times\textbf{H})\cdot B-\nabla\times(B\times D)\cdot\textbf{H}\right) dxdt\nonumber\\
&&+\int_{0}^{\tau}\int_{\Omega}\left(-(D\times B)\cdot(\nabla\times D)+(\textbf{H}\times\textbf{u})\cdot(\nabla\times D)\right) dxdt\nonumber\\
&=&\int_{0}^{\tau}\int_{\Omega}((\textbf{H}-D)^{\top}\nabla B(\textbf{H}-D)+\frac{1}{2}\nabla(|\textbf{H}-D|^2)\cdot B)dxdt\nonumber\\
&&+\int_{0}^{\tau}\int_{\Omega}(D-\textbf{H})^{\top}\nabla(\textbf{u}-B)D+\frac{1}{2}\nabla(D(D-\textbf{H}))(\textbf{u}-B)dxdt\nonumber\\
&&-\int_{0}^{\tau}\int_{\Omega}(D^{\top}\nabla(\textbf{u}-B)D+\frac{1}{2}\nabla(|D|^2)(\textbf{u}-B))dxdt.
\end{eqnarray}
Note that
\begin{eqnarray*}
\partial_y(\partial_{\rho}H_{C}(A,C))&=&-s(A,C)\partial_yC-A\partial_{\rho}s(A,C)\partial_yC+\partial^2_{\rho,\rho}H_C(A,C)\partial_y\rho\nonumber\\
&&+\partial^2_{\rho,\theta}H_C(A,C)\partial_yC,~~for~y=t,x.
\end{eqnarray*}
Thus it follows from (\ref{E2-R5}) and (\ref{E2-8}) that
\begin{eqnarray}\label{E2-6}
&&\int_{\Omega}(\frac{1}{2}\rho|\textbf{u}-B|^2+\frac{1}{2}|\textbf{H}-D|^2+H_{C}(\rho,\theta)-\partial_{\rho}(H_{C})(A,C)(\rho-A)-H_{C}(A,C))(\tau,\cdot)dx\nonumber\\
&&+\int_{0}^{\tau}\int_{\Omega}\left(|\nabla\textbf{u}|^2+\nu|\nabla\times D-\nabla\times\textbf{H}|^2\right)dxdt+\int_{0}^{\tau}\int_{\Omega}\frac{C}{\theta}(\textbf{S}(\theta,\nabla_x\textbf{u}):\nabla_x\textbf{u}-\frac{\textbf{q}(\theta,\nabla_x\theta)\cdot\nabla_x\theta}{\theta})dxdt\nonumber\\
&\leq&\int_{\Omega}(\frac{1}{2}\rho_0|\textbf{u}_0-B(0,\cdot)|^2+\frac{1}{2}|\textbf{H}_0-D_0|^2+(H_{C(0,\cdot)}(\rho(0,\cdot),\theta(0,\cdot))\nonumber\\
&&-\partial_{\rho}(H_{C(0,\cdot)})(A(0,\cdot),C(0,\cdot))(\rho_0-A(0,\cdot))-H_{C(0,\cdot)}(A(0,\cdot),C(0,\cdot))))dx\nonumber\\
&&+\int_{0}^{\tau}\int_{\Omega}((\rho\partial_tB+\rho\textbf{u}\cdot\nabla_xB)\cdot(B-\textbf{u})-P(\rho,\theta)\textbf{div}_xB+\textbf{S}(\theta,\nabla_x\textbf{u}):\nabla_xB)dxdt\nonumber\\
&&+\int_{0}^{\tau}\int_{\Omega}((\textbf{H}-D)^{\top}\nabla B(\textbf{H}-D)+\frac{1}{2}\nabla(|\textbf{H}-D|^2)\cdot B)dxdt\nonumber\\
&&+\int_{0}^{\tau}\int_{\Omega}(D-\textbf{H})^{\top}\nabla(\textbf{u}-B)D+\frac{1}{2}\nabla(D(D-\textbf{H}))(\textbf{u}-B)dxdt\nonumber\\
&&-\int_{0}^{\tau}\int_{\Omega}(D^{\top}\nabla(\textbf{u}-B)D+\frac{1}{2}\nabla(|D|^2)(\textbf{u}-B))dxdt\nonumber\\
&&-\int_{0}^{\tau}\int_{\Omega}(\rho(s(\rho,\theta)-s(A,C))\partial_tC+\rho(s(\rho,\theta)-s(A,C))\textbf{u}\cdot\nabla_x C+\frac{\textbf{q}(\theta,\nabla_x\theta)}{\theta}\cdot\nabla_x C)dxdt\nonumber\\
&&+\int_0^{\tau}\int_{\Omega}\rho(A\partial_{\rho}s(A,C)\partial_tC+r\partial_{\rho}s(A,C)\textbf{u}\cdot\nabla_xC)dxdt\nonumber\\
&&-\int_0^{\tau}\int_{\Omega}\rho(\partial^2_{\rho,\rho}(H_C)(A,C)\partial_tA+\partial^2_{\rho,\theta}(H_{C})(A,C)\partial_tC)dxdt\nonumber\\
&&-\int_0^{\tau}\int_{\Omega}\rho\textbf{u}(\partial^2_{\rho,\rho}(H_C)(A,C)\nabla_xA+\partial^2_{\rho,\theta}(H_{C})(A,C)\nabla_xC)dxdt\nonumber\\
&&+\int_0^{\tau}\int_{\Omega}\partial_t(A\partial_{\rho}(H_{C})(A,C)-H_{C}(A,C))dxdt.
\end{eqnarray}
Following \cite{Car,Fei1,Sa}, introducing the quantity as
\begin{eqnarray*}
\Gamma(\rho,\theta|C,C)=H_{C}(\rho,\theta)-\partial_{\rho}H_{C}(A,C)(\rho-A)-H_{C}(A,C),
\end{eqnarray*}
where
\begin{eqnarray*}
H_C(\rho,\theta)=\rho e(\rho,\theta)-C\rho s(\rho,\theta).
\end{eqnarray*}
Note that
\begin{eqnarray}\label{E2-12}
&&\partial^2_{\rho,\rho}H_C(A,C)=\frac{1}{A}\partial_{\rho}P(A,C),~~A\partial_{\rho}s(A,C)=-\frac{1}{C}\partial_{\theta}P(A,C),\nonumber\\
&&\partial^2_{\rho,\theta}H_C(A,C)=\partial_{\rho}(\rho(\theta-C)\partial_{\theta}s)(A,C)=(\theta-C)\partial_{\rho}(\rho\partial_{\theta}s(\rho,\theta))(A,C)=0,~~~~~~\\
&&A\partial_{\rho}(H_{C})(A,C)-H_{C}(A,C)=P(A,C).\nonumber
\end{eqnarray}
Therefore, we can obtain a kind of relative entropy inequality by simplifying (\ref{E2-6}) as
\begin{eqnarray}\label{E2-7}
&&\int_{\Omega}(\frac{1}{2}\rho|\textbf{u}-B|^2+\frac{1}{2}|\textbf{H}-D|^2+\Gamma(\rho,\theta|A,C))(\tau,\cdot)dx+\nu\int_{0}^{\tau}\int_{\Omega}|\nabla\times D-\nabla\times\textbf{H}|^2dxdt\nonumber\\
&&+\int_{0}^{\tau}\int_{\Omega}\frac{C}{\theta}(\textbf{S}(\theta,\nabla_x\textbf{u}):\nabla_x\textbf{u}-\frac{\textbf{q}(\theta,\nabla_x\theta)\cdot\nabla_x\theta}{\theta})dxdt\nonumber\\
&\leq&\int_{\Omega}(\frac{1}{2}\rho_0|\textbf{u}_0-B(0,\cdot)|^2+\frac{1}{2}|\textbf{H}_0-D_0(0,\cdot)|^2+\Gamma(\rho_0,\theta_0|A(0,\cdot),C(0,\cdot)))\nonumber\\
&&-\int_{0}^{\tau}\int_{\Omega}\rho(\textbf{u}-B)\cdot\nabla_xB\cdot(B-\textbf{u})dxdt+\int_{0}^{\tau}\int_{\Omega}\rho(s(\rho,\theta)-s(A,C))(B-\textbf{u})\cdot\nabla_xCdxdt\nonumber\\
&&+\int_{0}^{\tau}\int_{\Omega}((\rho\partial_tB+\rho\textbf{u}\cdot\nabla_xB)\cdot(B-\textbf{u})-P(\rho,\theta)\textbf{div}_xB+\textbf{S}(\theta,\nabla_x\textbf{u}):\nabla_xB)dxdt\nonumber\\
&&+\int_{0}^{\tau}\int_{\Omega}((\textbf{H}-D)^{\top}\nabla B(\textbf{H}-D)+\frac{1}{2}\nabla(|\textbf{H}-D|^2)\cdot B)dxdt\nonumber\\
&&+\int_{0}^{\tau}\int_{\Omega}(D-\textbf{H})^{\top}\nabla(\textbf{u}-B)D+\frac{1}{2}\nabla(D(D-\textbf{H}))(\textbf{u}-B)dxdt\nonumber\\
&&-\int_{0}^{\tau}\int_{\Omega}(D^{\top}\nabla(\textbf{u}-B)D+\frac{1}{2}\nabla(|D|^2)(\textbf{u}-B))dxdt\nonumber\\
&&-\int_{0}^{\tau}\int_{\Omega}(\rho(s(\rho,\theta)-s(A,C))\partial_tC+\rho(s(\rho,\theta)-s(A,C))\textbf{u}\cdot\nabla_x C+\frac{\textbf{q}(\theta,\nabla_x\theta)}{\theta}\cdot\nabla_x C)dxdt\nonumber\\
&&+\int_0^{\tau}\int_{\Omega}((1-\frac{\rho}{A})\partial_tP(A,C)-\frac{\rho}{A}\textbf{u}\cdot\nabla_xP(A,C))dxdt.
\end{eqnarray}
Now we state the weak-strong uniqueness property to the full magnetohydrodynamic system (\ref{E1-1})-(\ref{E1-1R})
on a bounded Lipschitz domains with Dirichlet
boundary conditions.
\begin{theorem}
Let $\Omega\subset\textbf{R}^3$ be a bounded Lipschitz domain and $(\rho,\textbf{u},\theta,\textbf{H})$ be a weak solution of the full magnetohydrodynamic system (\ref{E1-1})-(\ref{E1-1R}) in $(0,T)\times\Omega$ and $(\rho',\textbf{u}',\theta',\textbf{H}')$ be a strong solution emanating from the same initial data (\ref{E2-x2}).
Assume that the thermodynamic
functions $P,$ $e$, $s$ satisfy hypotheses (\ref{E1-15})-(\ref{E1-22}), and that the transport
coefficients $\mu$, $\eta$ and $\kappa$ satisfy (\ref{E1-23})-(\ref{E1-24}). Then
\begin{eqnarray*}
\rho\equiv\rho',~~\textbf{u}=\textbf{u}',~~\theta=\theta',~~\textbf{H}=\textbf{H}'.
\end{eqnarray*}
\end{theorem}

\section{Proof of Theorem 1}
In this section, we apply the relative entropy inequality to finish the proof of Theorem 1. Assume that $(\rho',\textbf{u}',\theta',\textbf{H}')$ is a classical (strong) solution to the full magnetohydrodynamic system in $(0,T)\times\Omega$, it satisfies that
\begin{eqnarray*}
\rho'(0,\cdot)=\rho_0,~~\textbf{u}'(0,\cdot)=\textbf{u}_0,~~\theta'(0,\cdot)=\theta_0,~~\textbf{H}'(0,\cdot)=\textbf{H}_0.
\end{eqnarray*}
Following \cite{Fei2,Fei1}, we introduce essential and residual component of
each quantity appearing in (\ref{E2-12}).
Thermodynamic stability hypothesis (\ref{E1-18R}) implies that
$\rho\mapsto H_{\theta'}(\rho,\theta')$ is strictly convex,
while $\theta\mapsto H_{\theta'}(\rho,\theta')$ attains its global minimum at $\theta=\theta'$.
Thus it has
\begin{eqnarray}\label{E3-2}
\Gamma(\rho,\theta|\rho',\theta')\geq c\left\{
\begin{array}{lll}
&&|\rho-\rho'|^2+|\theta-\theta'|^2~~if~~(\rho,\theta)\in[\rho_0',\rho_1']\times[\theta_0',\theta_1']\\
&&1+|\rho s(\rho,\theta)|+\rho e(\rho,\theta)~~otherwise,
\end{array}
\right.
\end{eqnarray}
where $[\rho',\theta']\in[\rho_0',\rho_1']\times[\theta_0',\theta_1']$, the constant $c$ depends on positive constants $\rho_0',\rho_1',\theta_0',\theta_1'$ and the structural properties of the thermodynamic function $e$, $s$. More precisely, the restriction of positive constants $\rho_0',\rho_1',\theta_0',\theta_1'$ can be found in \cite{Fei1}.

Thus we can write each measurable function $h=h_{ess}+h_{res}$,
where
\begin{eqnarray*}
h_{ess}=\left\{
\begin{array}{lll}
&&h(t,x)~~if~~(\rho,\theta)\in[\rho_0',\rho_1']\times[\theta_0',\theta_1']\\
&&0~~otherwise.
\end{array}
\right.
\end{eqnarray*}
Taking $(A,B,C,D)=(\rho',\textbf{u}',\theta',\textbf{H}')$ in (\ref{E2-7}). By the fact that the initial data coincide, we have
\begin{eqnarray}\label{E3-1}
&&\int_{\Omega}(\frac{1}{2}\rho|\textbf{u}-\textbf{u}'|^2+\frac{1}{2}|\textbf{H}-\textbf{H}'|^2+\Gamma(\rho,\theta|\rho',\theta))(\tau,\cdot)dx+\nu\int_{0}^{\tau}\int_{\Omega}|\nabla\times\textbf{H}'-\nabla\times\textbf{H}|^2dxdt\nonumber\\
&&+\int_{0}^{\tau}\int_{\Omega}\frac{\theta'}{\theta}(\textbf{S}(\theta,\nabla_x\textbf{u}):\nabla_x\textbf{u}-\frac{\textbf{q}(\theta,\nabla_x\theta)\cdot\nabla_x\theta}{\theta})dxdt\nonumber\\
&\leq&\int_{\Omega}(\frac{1}{2}\rho_0|\textbf{u}_0-\textbf{u}'(0,\cdot)|^2+\frac{1}{2}|\textbf{H}_0-\textbf{H}'_0(0,\cdot)|^2+\Gamma(\rho_0,\theta_0|\rho'(0,\cdot),\theta'(0,\cdot)))\nonumber\\
&&+\int_{0}^{\tau}\int_{\Omega}\rho|\textbf{u}-\textbf{u}'|^2|\nabla_x\textbf{u}'|dxdt+\int_{0}^{\tau}\int_{\Omega}\rho(s(\rho,\theta)-s(\rho,\theta))(\textbf{u}'-\textbf{u})\cdot\nabla_x\theta'dxdt\nonumber\\
&&+\int_{0}^{\tau}\int_{\Omega}((\rho\partial_t\textbf{u}'+\rho\textbf{u}\cdot\nabla_x\textbf{u}')\cdot(\textbf{u}'-\textbf{u})-P(\rho,\theta)\textbf{div}_x\textbf{u}'+\textbf{S}(\theta,\nabla_x\textbf{u}):\nabla_x\textbf{u}')dxdt\nonumber\\
&&+\int_{0}^{\tau}\int_{\Omega}((\textbf{H}-\textbf{H}')^{\top}\nabla\textbf{u}'(\textbf{H}-\textbf{H}')+\frac{1}{2}\nabla(|\textbf{H}-\textbf{H}'|^2)\cdot \textbf{u}')dx\nonumber\\
&&+\int_{0}^{\tau}\int_{\Omega}(\textbf{H}'-\textbf{H})^{\top}\nabla(\textbf{u}-\textbf{u}')\textbf{H}'+\frac{1}{2}\nabla(\textbf{H}'(\textbf{H}'-\textbf{H}))(\textbf{u}-\textbf{u}')dxdt\nonumber\\
&&-\int_{0}^{\tau}\int_{\Omega}(\textbf{H}'^{\top}\nabla(\textbf{u}-\textbf{u}')\textbf{H}'+\frac{1}{2}\nabla(|\textbf{H}'|^2)(\textbf{u}-\textbf{u}'))dxdt\nonumber\\
&&-\int_{0}^{\tau}\int_{\Omega}(\rho(s(\rho,\theta)-s(\rho',\theta'))\partial_t\theta'+\rho(s(\rho,\theta)-s(\rho',\theta'))\textbf{u}'\cdot\nabla_x\theta'+\frac{\textbf{q}(\theta,\nabla_x\theta)}{\theta}\cdot\nabla_x\theta')dxdt\nonumber\\
&&+\int_0^{\tau}\int_{\Omega}((1-\frac{\rho}{\rho'})\partial_tP(\rho',\theta')-\frac{\rho}{\rho'}\textbf{u}\cdot\nabla_xP(\rho',\theta'))dxdt.
\end{eqnarray}
In what follows, we estimate the right-hand side of (\ref{E3-1}). It is easy to see that
\begin{eqnarray}\label{E3-3}
\int_{\Omega}\rho|\textbf{u}-\textbf{u}'|^2|\nabla_x\textbf{u}'|dx\leq\|\nabla_x\textbf{u}'\|_{\textbf{L}^{\infty}(\Omega;\textbf{R}^3)}\int_{\Omega}\rho|\textbf{u}-\textbf{u}'|^2dx.
\end{eqnarray}
By virtue of (\ref{E3-2}), using interpolation inequality, for any $\epsilon>0$, we derive
\begin{eqnarray}\label{E3-4}
&&\int_{\Omega}\rho(s(\rho,\theta)-s(\rho,\theta))(\textbf{u}'-\textbf{u})\cdot\nabla_x\theta'dx\nonumber\\
&\leq&2\rho_1'\|\nabla_x\theta'\|_{\textbf{L}^{\infty}(\Omega;\textbf{R}^3)}(\epsilon\|\textbf{u}'-\textbf{u}\|^2_{\textbf{L}^2(\Omega;\textbf{R}^3)}+c(\epsilon)\int_{\Omega}\Gamma(\rho,\theta|\rho',\theta')dx)\nonumber\\
&&+\|\nabla_x\theta'\|_{\textbf{L}^{\infty}(\Omega;\textbf{R}^3)}(\epsilon\|\textbf{u}'-\textbf{u}\|^2_{\textbf{L}^6(\Omega;\textbf{R}^3)}+c(\epsilon)\|[\rho(s(\rho,\theta)-s(\rho',\theta'))]_{res}\|^2_{\textbf{L}^{\frac{6}{5}}(\Omega)}).~~~~~
\end{eqnarray}
It follows from (\ref{E1-20})-(\ref{E1-21}) that
\begin{eqnarray}\label{E3-5}
|[\rho(s(\rho,\theta)-s(\rho',\theta'))]_{res}|
\leq c(\rho+\rho[\log\theta]^++\rho|\log\rho|+\theta^3).
\end{eqnarray}
Using (\ref{E1-18}), (\ref{E1-19})-(\ref{E1-20}),
\begin{eqnarray}\label{E3-6}
\rho e(\rho,\theta)\geq c(\rho^{\frac{5}{3}}+\theta^4),
\end{eqnarray}
and (\ref{E2-3})-(\ref{E2-x1}) imply
\begin{eqnarray}\label{E3-7}
t\mapsto\int_{\Omega}\Gamma(\rho,\theta|\rho',\theta')dx\in\textbf{L}^{\infty}(0,T).
\end{eqnarray}
By (\ref{E3-2}), (\ref{E3-5})-(\ref{E3-7}) and H\"{o}lder inequality,
\begin{eqnarray*}
\|[\rho(s(\rho,\theta)-s(\rho',\theta'))]_{res}\|^2_{\textbf{L}^{\frac{6}{5}}(\Omega)}\leq c(\int_{\Omega}\Gamma(\rho,\theta|\rho',\theta,))^{\frac{5}{3}}.
\end{eqnarray*}
So by (\ref{E3-4}), for any $\epsilon>0$, we obtain
\begin{eqnarray}\label{E3-8}
\int_{\Omega}\rho(s(\rho,\theta)-s(\rho,\theta))(\textbf{u}'-\textbf{u})\cdot\nabla_x\theta'dx\leq\epsilon\|\textbf{u}'-\textbf{u}\|^2_{\textbf{W}_0^{1,2}(\Omega;\textbf{R}^3)}+c'(\epsilon,\cdot)\int_{\Omega}\Gamma(\rho,\theta|\rho',\theta')dx,~~~
\end{eqnarray}
where $c'(\epsilon,\cdot)$ is a generic constant depending on $\epsilon$, $\rho'$, $\textbf{u}'$ and $\theta'$ through the norms induced by (\ref{E2-x1})-(\ref{E2-x2}), while $c'(\cdot)$ is independent of $\epsilon$ but depends on $\rho'$, $\textbf{u}'$, $\theta'$, $\rho_0'$ and $\theta_0,$ through the norms induced by (\ref{E2-x1})-(\ref{E2-x2}).

Similar with estimating (\ref{E3-8}), we get
\begin{eqnarray}\label{E3-9}
&&\int_{\Omega}\frac{1}{\rho'}(\rho-\rho')(\textbf{u}'-\textbf{u}')\cdot(\textbf{div}_x\textbf{S}(\theta',\nabla_x\textbf{u})-\nabla_xP(\rho',\theta')+(\nabla\times\textbf{H}')\times\textbf{H}'))dx\nonumber\\
&=&\int_{\Omega}[\rho'^{-1}(\rho-\rho')(\textbf{u}'-\textbf{u}')\cdot[\textbf{div}_x\textbf{S}(\theta',\nabla_x\textbf{u})-\nabla_xP(\rho',\theta')+(\nabla\times\textbf{H}')\times\textbf{H}')]_{ess}dx\nonumber\\
&&+\int_{\Omega}[\rho'^{-1}(\rho-\rho')(\textbf{u}'-\textbf{u}')\cdot[\textbf{div}_x\textbf{S}(\theta',\nabla_x\textbf{u})-\nabla_xP(\rho',\theta')+(\nabla\times\textbf{H}')\times\textbf{H}')]_{ess}dx\nonumber\\
&\leq&c'(\epsilon,\cdot)\|[\rho-\rho']_{ess}\|^2_{\textbf{L}^2(\Omega)}+\epsilon\|\textbf{u}'-\textbf{u}\|^2_{\textbf{L}^2(\Omega;\textbf{R}^3)}\nonumber\\
&&+c'(\epsilon,\cdot)(\|[\rho]_{ess}\|^2_{\textbf{L}^{\frac{6}{5}}(\Omega)}+\|[1]_{res}\|^2_{\textbf{L}^{\frac{6}{5}}(\Omega)})
+\epsilon\|\textbf{u}'-\textbf{u}\|^2_{\textbf{L}^6(\Omega;\textbf{R}^3)}.
\end{eqnarray}
So using integrating by parts, then virtue of (\ref{E3-2}), (\ref{E3-7}), (\ref{E3-9}) and $\textbf{W}^{1,2}(\Omega)\hookrightarrow\textbf{L}^6(\Omega)$, we derive
\begin{eqnarray}\label{E3-10}
&&|\int_{\Omega}\rho(\partial_t\textbf{u}'+\textbf{u}'\cdot\nabla_x\textbf{u}')\cdot(\textbf{u}'-\textbf{u})dx|\nonumber\\
&=&|\int_{\Omega}\frac{\rho}{\rho'}(\textbf{u}'-\textbf{u})\cdot(\textbf{div}_x\textbf{S}(\theta',\nabla_x\textbf{u}')-\nabla_xP(\rho',\theta')+(\nabla\times\textbf{H}')\times\textbf{H}')dx|\nonumber\\
&\leq&\int_{\Omega}|\frac{\rho-\rho'}{\rho'}(\textbf{u}'-\textbf{u})\cdot(\textbf{div}_x\textbf{S}(\theta',\nabla_x\textbf{u}')-\nabla_xP(\rho',\theta')+(\nabla\times\textbf{H}')\times\textbf{H}')|dx\nonumber\\
&&+|\int_{\Omega}(\textbf{u}'-\textbf{u})\cdot(\textbf{div}_x\textbf{S}(\theta',\nabla_x\textbf{u}')-\nabla_xP(\rho',\theta')+(\nabla\times\textbf{H}')\times\textbf{H}')dx|\nonumber\\
&\leq&c'(\epsilon,\cdot)\|[\rho-\rho']_{ess}\|^2_{\textbf{L}^2(\Omega)}+\epsilon\|\textbf{u}'-\textbf{u}\|^2_{\textbf{L}^2(\Omega;\textbf{R}^3)}+c'(\epsilon,\cdot)(\|[\rho]_{ess}\|^2_{\textbf{L}^{\frac{6}{5}}(\Omega)}+\|[1]_{res}\|^2_{\textbf{L}^{\frac{6}{5}}(\Omega)})\nonumber\\
&&+\epsilon\|\textbf{u}'-\textbf{u}\|^2_{\textbf{L}^6(\Omega;\textbf{R}^3)}\nonumber\\
&&+|\int_{\Omega}\left(\textbf{S}(\theta',\nabla_x\textbf{u}'):\nabla_x(\textbf{u}'-\textbf{u})+P(\rho',\theta')\textbf{div}_x(\textbf{u}'-\textbf{u})+((\nabla\times\textbf{H}')\times\textbf{H}')\cdot(\textbf{u}'-\textbf{u})\right)dx|\nonumber\\
&\leq&|\int_{\Omega}\left(\textbf{S}(\theta',\nabla_x\textbf{u}'):\nabla_x(\textbf{u}'-\textbf{u})+P(\rho',\theta')\textbf{div}_x(\textbf{u}'-\textbf{u})+((\nabla\times\textbf{H}')\times\textbf{H}')\cdot(\textbf{u}'-\textbf{u})\right)dx|\nonumber\\
&&+\left(\epsilon\|\textbf{u}-\textbf{u}'\|_{\textbf{W}_0^{1,2}(\Omega;\textbf{R}^3)}+c(\epsilon)\int_{\Omega}\Gamma(\rho,\theta|\rho',\theta')dx\right).
\end{eqnarray}
By H\"{o}lder inequality and (\ref{E1-11})-(\ref{E1-12}), we derive
\begin{eqnarray}\label{E9-1}
\int_{0}^{\tau}\int_{\Omega}((\textbf{H}-\textbf{H}')^{\top}\nabla\textbf{u}'(\textbf{H}-\textbf{H}')dxdt
\leq\int_{0}^{\tau}(\|\nabla\textbf{u}'\|_{\textbf{L}^{\infty}(\Omega;\textbf{R}^3)}\int_{\Omega}|\textbf{H}-\textbf{H}'|^2dx)dt,~~~~~~
\end{eqnarray}
\begin{eqnarray}\label{E9-2}
\frac{1}{2}\int_{0}^{\tau}\int_{\Omega}\nabla(|\textbf{H}-\textbf{H}'|^2)\cdot\textbf{u}'dxdxdt
\leq\frac{1}{2}\int_{0}^{\tau}(\|\nabla\textbf{u}'\|_{\textbf{L}^{\infty}(\Omega;\textbf{R}^3)}\int_{\Omega}|\textbf{H}-\textbf{H}'|^2dx)dt,~~~~~
\end{eqnarray}
\begin{eqnarray}\label{E9-3}
\int_{0}^{\tau}\int_{\Omega}(\textbf{H}'-\textbf{H})^{\top}\nabla(\textbf{u}-\textbf{u}')\textbf{H}'dxdt
&\leq&\int_{0}^{\tau}c_{\epsilon}\|\textbf{H}'\|_{\textbf{L}^{\infty}(\Omega;\textbf{R}^3)}
\int_{\Omega}|\textbf{H}'-\textbf{H}|^2dxdt\nonumber\\
&&+\int_0^{\tau}\epsilon\|\textbf{H}'\|_{\textbf{L}^{\infty}(\Omega;\textbf{R}^3)}\|\textbf{u}-\textbf{u}'\|^2_{\textbf{W}^{1,2}(\Omega;\textbf{R}^3)}dt,~~~~~~~
\end{eqnarray}
\begin{eqnarray}\label{E9-4}
\frac{1}{2}\int_{0}^{\tau}\int_{\Omega}\nabla(\textbf{H}'(\textbf{H}'-\textbf{H}))(\textbf{u}-\textbf{u}')dxdt
&\leq&\frac{1}{4}\int_{0}^{\tau}c_{\epsilon}\int_{\Omega}|\textbf{H}'-\textbf{H}|^2dxdt\nonumber\\
&&+\frac{1}{4}\int_0^{\tau}\epsilon\|\textbf{H}'\|^2_{\textbf{L}^{\infty}(\Omega;\textbf{R}^3)}\|\textbf{u}-\textbf{u}'\|^2_{\textbf{W}^{1,2}(\Omega;\textbf{R}^3)}dt,~~~~~
\end{eqnarray}
Thus combing with (\ref{E9-1})-(\ref{E9-4}) and (\ref{E9-5}), we have
\begin{eqnarray}\label{E3-11}
&&\int_{0}^{\tau}\int_{\Omega}((\textbf{H}-\textbf{H}')^{\top}\nabla\textbf{u}'(\textbf{H}-\textbf{H}')+\frac{1}{2}\nabla(|\textbf{H}-\textbf{H}'|^2)\cdot \textbf{u}')dx\nonumber\\
&&+\int_{0}^{\tau}\int_{\Omega}(\textbf{H}'-\textbf{H})^{\top}\nabla(\textbf{u}-\textbf{u}')\textbf{H}'+\frac{1}{2}\nabla(\textbf{H}'(\textbf{H}'-\textbf{H}))(\textbf{u}-\textbf{u}')dxdt\nonumber\\
&&-\int_{0}^{\tau}\int_{\Omega}(\textbf{H}'^{\top}\nabla(\textbf{u}-\textbf{u}')\textbf{H}'+\frac{1}{2}\nabla(|\textbf{H}'|^2)(\textbf{u}-\textbf{u}'))dxdt\nonumber\\
&\leq&\int_{0}^{\tau}\left(\epsilon c''\|\textbf{u}-\textbf{u}'\|^2_{\textbf{W}^{1,2}(\Omega;\textbf{R}^3)}+c_{\epsilon}\int_{\Omega}|\textbf{H}'-\textbf{H}|^2dx\right)dt\nonumber\\
&&+\int_0^{\tau}\int_{\Omega}((\nabla\times\textbf{H}')\times\textbf{H}')(\textbf{u}-\textbf{u}')dxdt,
\end{eqnarray}
where $c''=c''(\|\textbf{H}'\|^2_{\textbf{L}^{\infty}(\Omega;\textbf{R}^3)})$ and
$c_{\epsilon}=c_{\epsilon}(\|\nabla\textbf{u}'\|_{\textbf{L}^{\infty}(\Omega;\textbf{R}^3)},\|\textbf{H}'\|_{\textbf{L}^{\infty}(\Omega;\textbf{R}^3)})$ denote constants, $\epsilon>0$ sufficient small.

In what follows, we estimate the next term. Using Taylor-Lagrange formula, we derive
\begin{eqnarray}\label{E3-12}
&&\int_{\Omega}\rho(s(\rho,\theta)-s(\rho',\theta'))\partial_t\theta'dx\nonumber\\
&\leq&\int_{\Omega}\rho'[\partial_{\rho}s(\rho',\theta')(\rho-\rho')+\partial_{\theta}s(\rho',\theta')(\theta-\theta')]\partial_t\theta'dx+4c(\cdot)\int_{\Omega}\Gamma(\rho,\theta|\rho',\theta')dx.~~~~~~
\end{eqnarray}
Similar with getting (\ref{E3-12}) and (\ref{E3-4}), respectively, we have
\begin{eqnarray}\label{E3-13}
&&-\int_{\Omega}\rho(s(\rho,\theta)-s(\rho',\theta'))\textbf{u}'\cdot\nabla_x\theta'dx\nonumber\\
&\leq&-\int_{\Omega}\rho'[\partial_{\rho}s(\rho',\theta')(\rho-\rho')+\partial_{\theta}s(\rho',\theta')(\theta-\theta')]\textbf{u}'\cdot\nabla_x\theta'dx+c(\cdot)\int_{\Omega}\Gamma(\rho,\theta|\rho',\theta')dx~~~~~~
\end{eqnarray}
and
\begin{eqnarray}\label{E3-14}
&&|\int_{\Omega}\frac{\rho-\rho'}{\rho'}\nabla_xP(\rho',\theta')\cdot(\textbf{u}-\textbf{u}')dx|\nonumber\\
&\leq&c(|\nabla_x\rho',|\nabla_x\theta'|)\left(\epsilon\|\textbf{u}-\textbf{u}'\|^2_{\textbf{W}^{1,2}(\Omega;\textbf{R}^3)}+\int_{\Omega}\Gamma(\rho,\theta|\rho',\theta,)dx\right).~~~~~
\end{eqnarray}
Now we estimate the last term of right-hand side of (\ref{E3-1}). By (\ref{E3-14}),
\begin{eqnarray}\label{E3-15}
&&|\int_{\Omega}\left(\frac{\rho-\rho'}{\rho'}\partial_tP(\rho',\theta')-\frac{\rho}{\rho'}\textbf{u}\cdot\nabla_xP(\rho',\theta')\right)dx\nonumber\\
&\leq&\int_{\Omega}\frac{\rho-\rho'}{\rho'}(\partial_tP(\rho',\theta')+\textbf{u}\cdot\nabla_xP(\rho',\theta'))dx
+\int_{\Omega}P(\rho',\theta')\textbf{div}_x\textbf{u}dx\nonumber\\
&&+c(|\nabla_x\rho',|\nabla_x\theta'|)\left(\epsilon\|\textbf{u}-\textbf{u}'\|^2_{\textbf{W}^{1,2}(\Omega;\textbf{R}^3)}+\int_{\Omega}\Gamma(\rho,\theta|\rho',\theta,)dx\right).~~
\end{eqnarray}
Thus by (\ref{E3-1})-(\ref{E3-3}), (\ref{E3-10})-(\ref{E3-12}) and (\ref{E3-15}), we obtain the following relative entropy inequality
\begin{eqnarray}\label{E3-16}
&&\int_{\Omega}(\frac{1}{2}\rho|\textbf{u}-\textbf{u}'|^2+\frac{1}{2}|\textbf{H}-\textbf{H}'|^2+\Gamma(\rho,\theta|\rho',\theta))(\tau,\cdot)dx+\nu\int_{0}^{\tau}\int_{\Omega}|\nabla\times\textbf{H}'-\nabla\times\textbf{H}|^2dxdt\nonumber\\
&&+\int_{0}^{\tau}\int_{\Omega}(\frac{\theta'}{\theta}\textbf{S}(\theta,\nabla_x\textbf{u}):\nabla_x\textbf{u}
-\textbf{S}(\theta',\nabla_x\textbf{u}'):(\nabla_x\textbf{u}-\nabla_x\textbf{u}')-\textbf{S}(\theta,\nabla_x\textbf{u}):\nabla_x\textbf{u}'dxdt)\nonumber\\
&&+\int_{0}^{\tau}\int_{\Omega}\left(\frac{\textbf{q}(\theta,\nabla_x\theta)\cdot\nabla_x\theta'}{\theta}-\frac{\theta'}{\theta}\frac{\textbf{q}(\theta,\nabla_x\theta)\cdot\nabla_x\theta}{\theta}dxdt \right)\nonumber\\
&\leq&\int_{0}^{\tau}\left(\epsilon\|\textbf{u}'-\textbf{u}\|^2_{\textbf{W}_0^{1,2}(\Omega;\textbf{R}^3)}+c'(\epsilon,\cdot)\int_{\Omega}(\Gamma(\rho,\theta|\rho',\theta')+\frac{\rho}{2}|\textbf{u}-\textbf{u}'|^2)dx\right)dt\nonumber\\
&&+\int_0^{\tau}\int_{\Omega}(P(\rho',\theta')-P(\rho,\theta))\textbf{div}_x\textbf{u}'+((\nabla\times\textbf{H}')\times\textbf{H}')\cdot(\textbf{u}'-\textbf{u})dxdt\nonumber\\
&&+\int_{0}^{\tau}\left(\epsilon c''\|\textbf{u}-\textbf{u}'\|^2_{\textbf{W}^{1,2}(\Omega;\textbf{R}^3)}+c_{\epsilon}\int_{\Omega}|\textbf{H}'-\textbf{H}|^2dx\right)dt\nonumber\\
&&-\int_0^{\tau}\int_{\Omega}\rho'[\partial_{\rho}s(\rho',\theta')(\rho-\rho')+\partial_{\theta}s(\rho',\theta')(\theta-\theta')](\partial_t\theta'+\textbf{u}'\cdot\nabla_x\theta')dxdt\nonumber\\
&&+\int_0^{\tau}\int_{\Omega}\frac{\rho-\rho'}{\rho'}(\partial_tP(\rho',\theta')+\textbf{u}\cdot\nabla_xP(\rho',\theta'))dxdt.
\end{eqnarray}
Next we simplify the relative entropy inequality (\ref{E3-16}). Using (\ref{E2-12}) and the fact that $\rho'$ and $\textbf{u}'$ satisfy the equation of continuity (\ref{E1-1}), we derive
\begin{eqnarray}\label{E3-17}
&&-\int_{\Omega}\rho'[\partial_{\rho}s(\rho',\theta')(\rho-\rho')+\partial_{\theta}s(\rho',\theta')(\theta-\theta')](\partial_t\theta'+\textbf{u}'\cdot\nabla_x\theta')dx+\int_{\Omega}\frac{\rho-\rho'}{\rho'}(\partial_tP(\rho',\theta')+\textbf{u}\cdot\nabla_xP(\rho',\theta'))dx\nonumber\\
&=&-\int_{\Omega}\rho'(\theta-\theta')\left(\frac{1}{\theta}(\textbf{S}(\theta',\nabla_x\textbf{u}'):\nabla_x\textbf{u}'-\frac{\textbf{q}(\theta',\nabla_x\theta')\cdot\nabla_x\theta'}{\theta'})-\textbf{div}_x(\frac{\textbf{q}(\theta',\nabla_x\theta')}{\theta})\right)dx\nonumber\\
&&+\int_{\Omega}\left((\theta-\theta')\partial_{\theta}P(\rho',\theta')+(\rho-\rho')\partial_{\rho}P(\rho',\theta')\right)\textbf{div}_x\textbf{u}'dx.
\end{eqnarray}
Note that
\begin{eqnarray}\label{E3-18R}
&&|\int_{\Omega}\left(P(\rho',\theta')-\partial_{\rho}P(\rho',\theta')(\rho'-\rho)-\partial_{\theta}P(\rho',\theta')(\theta'-\theta)-P(\rho,\theta)\textbf{div}_x\textbf{u}'\right)dx|\nonumber\\
&\leq&c\|\textbf{div}_x\textbf{u}\|_{\textbf{L}^{\infty}(\Omega)}\int_{\Omega}\Gamma(\rho,\theta|\rho',\theta')dx.
\end{eqnarray}
Thus by (\ref{E3-16})-(\ref{E3-17}) and (\ref{E3-18R}), for any $\epsilon>0$, we obtain
\begin{eqnarray}\label{E3-18}
&&\int_{\Omega}(\frac{1}{2}\rho|\textbf{u}-\textbf{u}'|^2+\frac{1}{2}|\textbf{H}-\textbf{H}'|^2+\Gamma(\rho,\theta|\rho',\theta))(\tau,\cdot)dx+\nu\int_{0}^{\tau}\int_{\Omega}|\nabla\times\textbf{H}'-\nabla\times\textbf{H}|^2dxdt\nonumber\\
&&+\int_{0}^{\tau}\int_{\Omega}(\frac{\theta'}{\theta}\textbf{S}(\theta,\nabla_x\textbf{u}):\nabla_x\textbf{u}
-\textbf{S}(\theta',\nabla_x\textbf{u}'):(\nabla_x\textbf{u}-\nabla_x\textbf{u}')-\textbf{S}(\theta,\nabla_x\textbf{u}):\nabla_x\textbf{u}'+\frac{\theta-\theta'}{\theta'}\textbf{S}(\theta',\nabla_x\textbf{u}':\nabla_x\textbf{u}')dxdt)\nonumber\\
&&+\int_{0}^{\tau}\int_{\Omega}(\frac{\textbf{q}(\theta,\nabla_x\theta)\cdot\nabla_x\theta'}{\theta}-\frac{\theta'}{\theta}\frac{\textbf{q}(\theta,\nabla_x\theta)\cdot\nabla_x\theta}{\theta}+(\theta'-\theta)\frac{\textbf{q}(\theta',\nabla_x\theta')}{\theta'^2}+\frac{\textbf{q}(\theta',\nabla_x\theta')}{\theta'}\cdot\nabla_x(\theta-\theta'))dxdt\nonumber\\
&\leq&\int_{0}^{\tau}\left(\epsilon\|\textbf{u}'-\textbf{u}\|^2_{\textbf{W}_0^{1,2}(\Omega;\textbf{R}^3)}+c'(\epsilon,\cdot)\int_{\Omega}(\Gamma(\rho,\theta|\rho',\theta')+\frac{1}{2}\rho|\textbf{u}-\textbf{u}'|^2)dx\right)dt\nonumber\\
&&+\int_{0}^{\tau}\left(\epsilon c''\|\textbf{u}-\textbf{u}'\|^2_{\textbf{W}^{1,2}(\Omega;\textbf{R}^3)}+c_{\epsilon}\int_{\Omega}|\textbf{H}'-\textbf{H}|^2dx\right)dt\nonumber\\
&&+\int_{\Omega}\left(P(\rho',\theta')-\partial_{\rho}P(\rho',\theta')(\rho'-\rho)-\partial_{\theta}P(\rho',\theta')(\theta'-\theta)-P(\rho,\theta)\textbf{div}_x\textbf{u}'\right)dx\nonumber\\
&\leq&\int_{0}^{\tau}\left(\epsilon\|\textbf{u}'-\textbf{u}\|^2_{\textbf{W}_0^{1,2}(\Omega;\textbf{R}^3)}+c'(\epsilon,\cdot)\int_{\Omega}(\Gamma(\rho,\theta|\rho',\theta')+\frac{1}{2}\rho|\textbf{u}-\textbf{u}'|^2)dx\right)dt\nonumber\\
&&+\int_{0}^{\tau}\left(\epsilon c''\|\textbf{u}-\textbf{u}'\|^2_{\textbf{W}^{1,2}(\Omega;\textbf{R}^3)}+c_{\epsilon}\int_{\Omega}|\textbf{H}'-\textbf{H}|^2dx\right)dt.
\end{eqnarray}
The estimate of the terms on $\textbf{S}$ and $\textbf{q}$ has been founded in \cite{Fei1}. So we omit it.
To conclude that we get
\begin{eqnarray*}
&&\int_{\Omega}(\frac{1}{2}\rho|\textbf{u}-\textbf{u}'|^2+\frac{1}{2}|\textbf{H}-\textbf{H}'|^2+\Gamma(\rho,\theta|\rho',\theta))(\tau,\cdot)dx+\nu\int_{0}^{\tau}\int_{\Omega}|\nabla\times\textbf{H}'-\nabla\times\textbf{H}|^2dxdt\nonumber\\
&&+c_1\int_{0}^{\tau}\int_{\Omega}\left(|\nabla_x\textbf{u}'-\nabla_x\textbf{u}|^2+|\nabla_x\theta'-\nabla_x\theta|^2+|\nabla_x\log\theta'-\nabla_x\log\theta|^2\right)dxdt\nonumber\\
&\leq&c_2\int_{0}^{\tau}\int_{\Omega}\left(\Gamma(\rho,\theta|\rho',\theta')+\frac{1}{2}\rho|\textbf{u}-\textbf{u}'|^2\right)dxdt\nonumber\\
&&+\int_{0}^{\tau}\left(\epsilon c''\|\textbf{u}-\textbf{u}'\|^2_{\textbf{W}^{1,2}(\Omega;\textbf{R}^3)}+c_{\epsilon}\int_{\Omega}|\textbf{H}'-\textbf{H}|^2dx\right)dt,~~for~almost~all~\tau\in(0,T).
\end{eqnarray*}
Furthermore, for sufficient small $\epsilon>0$, we obtain
\begin{eqnarray*}
&&\int_{\Omega}(\frac{1}{2}\rho|\textbf{u}-\textbf{u}'|^2+\frac{1}{2}|\textbf{H}-\textbf{H}'|^2+\Gamma(\rho,\theta|\rho',\theta))(\tau,\cdot)dx+\nu\int_{0}^{\tau}\int_{\Omega}|\nabla\times\textbf{H}'-\nabla\times\textbf{H}|^2dxdt\nonumber\\
&&+(c_1-c''\epsilon)\int_{0}^{\tau}\int_{\Omega}\left(|\nabla_x\textbf{u}'-\nabla_x\textbf{u}|^2+|\nabla_x\theta'-\nabla_x\theta|^2+|\nabla_x\log\theta'-\nabla_x\log\theta|^2\right)dxdt\nonumber\\
&\leq&c_2\int_{0}^{\tau}\left(\int_{\Omega}(\Gamma(\rho,\theta|\rho',\theta')+\frac{1}{2}\rho|\textbf{u}-\textbf{u}'|^2+\frac{1}{2}|\textbf{H}'-\textbf{H}|^2)\right)dx
\end{eqnarray*}
which implies that
\begin{eqnarray*}
\rho\equiv\rho',~~\textbf{u}\equiv\textbf{u}',~~\theta\equiv\theta',~~\textbf{H}\equiv\textbf{H}'.
\end{eqnarray*}
This completes the proof.

\begin{acknowledgements}
The author also expresses his sincere thanks to the
anonymous referees for very careful reading and for providing many valuable
comments and suggestions which led to improvement of this paper.
The author expresses his sincerely thanks to prof Zhifei Zhang and prof Yong Li for their comments and suggestions.
The author is supported by NSFC No 11201172, 11071101, the post doctor fund 2012M510243, SRFDP Grant No 20120061120002 and the 985 Project of
Jilin University.
\end{acknowledgements}




\end{document}